\documentclass[10pt]{article} 
\usepackage{amsmath}
\usepackage{amssymb} 
\usepackage{amscd}
\usepackage{theorem} 
\usepackage{xspace}

\def\dim{{\mbox{dim}}}

\def\calp{{\cal P}} 
\def\cala{{\cal A}} 
\def\calb{{\cal B}}

\def\bbbone{\mbox{\rm 1\hspace {-.6em} l}}

\newtheorem{proposition}{PROPOSITION}
\newtheorem{corollary}{COROLLARY}
\begin{document}

\baselineskip=0.5cm

\begin{center} 
 \thispagestyle{empty}
{\large\bf HOMOGENEOUS ALGEBRAS, STATISTICS AND COMBINATORICS}
\end{center} 
\vspace{0.5cm}

\begin{center} 
 Michel DUBOIS-VIOLETTE
\footnote{Laboratoire de Physique Th\'eorique, UMR 8627, Universit\'e Paris XI,
B\^atiment 210, F-91 405 Orsay Cedex, France\\
Michel.Dubois-Violette$@$th.u-psud.fr\\
} and 
Todor POPOV
\footnote{Theoretical Physics Division, 
Institute for Nuclear Research and Nuclear Energy\\
Tsarigradsko Chaussee 72,
BG-1784, Sofia, Bulgarie\\
tpopov@inrne.bas.bg\\
{\sl and} Laboratoire de Physique Th\'eorique, UMR 8627, Universit\'e Paris XI,
B\^atiment 210, F-91 405 Orsay Cedex, France
}

\end{center} \vspace{1cm}

\begin{center} \today \end{center}

\vspace {0,5cm}

\begin{abstract}
After some generalities on homogeneous algebras, we give a formula connecting the Poincar\'e series of a homogeneous algebra with the homology of the corresponding Koszul complex generalizing thereby a standard result for quadratic algebras. We then investigate two particular types of cubic algebras: The first one called the parafermionic (parabosonic) algebra is the algebra generated by the creation operators of the universal fermionic (bosonic) parastatics with $D$ degrees of freedom while the second is the plactic algebra that is the algebra of the plactic monoid with entries in $\{1,2,\dots, D\}$. In the case $D=2$ we describe the relations with the cubic  Artin-Schelter algebras. It is pointed out that the natural action of $GL(2)$ on the parafermionic algebra for $D=2$ extends as an action of the quantum group $GL_{p,q}(2)$ on the generic cubic Artin-Schelter regular algebra of type $S_1$; $p$ and $q$ being related to the Artin-Schelter parameters. It is claimed that this has a counterpart for any integer $D\geq 2$.

\end{abstract}
\vspace{0,5cm}
{\bf MSC (2000)} : 58B34, 81R60, 16E65, 08C99, 05A17, 05E10.\\
{\bf Keywords} : Homogeneous algebras, duality, N-complexes, Koszul algebras, Young tableaux.

  \vspace{1cm}
  
\noindent LPT-ORSAY 02-60

\newpage

\section{Introduction}
 Parastatistics has been introduced  as a generalization of the canonical relations corresponding to Bose-Fermi alternative. It has a long history reviewed in \cite{FF} where many references can be found.
 Instead of the canonical relations
\begin{equation}
[a_k,a^\ast_\ell]_\pm = \delta_{k\ell} \bbbone
\label{eq1.1}
\end{equation}
\begin{equation}
[a_k,a_\ell]_\pm =0
\label{eq1.2}
\end{equation}
one takes the following relations of degree 3
\begin{equation}
[a_k,[a^\ast_\ell,a_m]_\pm]_-=2\delta_{k\ell} a_m
\label{eq1.3}
\end{equation}
\begin{equation}
[a_k,[a_\ell,a_m]_\pm]_-=0
\label{eq1.4}
\end{equation}
(which follow from (\ref{eq1.1}) and (\ref{eq1.2})) as basic relations.\\

Relations (\ref{eq1.3}) and (\ref{eq1.4}) are the parastatistics relations as described in \cite{Green}. The upper sign stands for (para) Bose statistics while the lower one stands for (para) Fermi statistics. From now on we concentrate on parafermionic systems with a finite number $D$ of degrees of freedom which means that in the above equations (\ref{eq1.3}) and (\ref{eq1.4}) we take the lower sign and that $k,\ell,m$ belong to the set $\{1,\dots,D\}$. Let $\calb$ denote the Borel subalgebra generated by the $a_k$. The relations (\ref{eq1.4}) read in our case 
\begin{equation}
[a_k,[a_\ell,a_m]]=0
\label{eq1.5}
\end{equation}
for any $k,\ell,m \in \{1,\dots,D\}$. By replacement of the $a_k$ by their adjoints $a^\ast_k$ (``creation operators") one can interpret $\calb$ as (spanning) the universal Fock space for fermionic parastatistics with $D$ degrees of freedom, which is the Fock space encompassing parafermions of arbitrary order. When no confusion arises $\calb$ will be refered to as {\sl the parafermionic algebra}.\\

A very useful tool in combinatorics is the structure of associative monoid which can be attached to the set of Young tableaux with entries in $\{1,\dots,D\}$ \cite{LS}, \cite{Ful}. The monoid so obtained is called {\sl the plactic monoid}. The algebra of the plactic monoid which will be refered to as {\sl the plactic algebra} is the algebra $\calp$ generated by $D$ elements $e_k$ ($k\in \{1,\dots,D\}$) with relations
\begin{equation}
\left.
\begin{array}{lllll}
e_\ell e_m e_k & = & e_\ell e_k e_m &  \mbox{if} &  k<\ell\leq m\\
\\
e_k e_m e_\ell & = & e_m e_k e_\ell &  \mbox{if} & k\leq \ell < m
\end{array}
\right\}
\label{eq1.6}
\end{equation}
for $k,\ell, m\in \{ 1,\dots, D\}$. These relations are the {\sl Knuth relations} \cite{Ful}.\\

Both algebras $\calb$ and $\calp$ are generated by $D$ elements with relations of degree 3, relations (\ref{eq1.5}) for $\calb$ and relations (\ref{eq1.6}) for $\calp$ : They are homogeneous algebras of degree 3 or cubic algebras \cite{BD-VW}. By giving  the degree 1 to their generators $\calb$ and $\calp$ are ($\mathbb N$-)graded algebras which furthermore admit both a homogeneous basis labelled by the Young tableaux with entries in $\{1,\dots,D\}$. This is clear by definition for $\calp$ and for $\calb$ this follows from the decomposition of the action by automorphisms of the linear group $GL(D)$ into irreducible components as will be explained in Section 3. It follows in particular that $\calb$ and $\calp$ have the same Poincar\'e series
\begin{equation}
P_\calb (t) = P_\calp(t)
\label{eq1.7}
\end{equation}
and share many common properties.\\

Our aim in the following is to analyse the algebras $\calb$ and $\calp$ by using the general tools developed in \cite{BD-VW} for homogeneous algebras which extend the classical ones for quadratic 
algebras \cite{Pri}, \cite{YuM2}. It is worth noticing here that homogeneous algebras constitute a class of algebras which is fundamental for the noncommutative version of algebraic geometry (see e.g. \cite{AS}, \cite{ATVB}, \cite{SmSt}, \cite{Tvdb}, \cite{RB3}, \cite{RB4}, \cite{connes:02}, \cite{connes:03}, \cite{AC.MDV}) and in which enters canonically the theory of $N$-complexes which itself received recent developments (see e.g. in \cite{BD-VW}, \cite{D-V2}, \cite{D-V4}, \cite{D-VH2}, \cite{D-VK}, \cite{D-VT2}, \cite{Kap}, \cite{KW}, \cite{Wam3}).  It should be stressed that this paper is not a paper of (para)statistics or of combinatorics but is a paper on homogeneous algebras in which the concepts and technics of \cite{BD-VW} are examplified by the analysis of the algebras $\calb$ and $\calp$ coming from these domains and it is expected that this will prove useful there.\\

Before this we shall complete the general analysis of \cite{BD-VW} by establishing a formula for the Poincar\'e series which generalizes to $N$-homogeneous algebras a result known for quadratic algebras and which is particularly simple and useful for Koszul algebras.\\

About the notations, we use throughout the notations of \cite{BD-VW} which are partly reviewed in Section 2. By an algebra we always mean here an associative unital algebra and when we speak of an algebra generated by some elements, the unit is not supposed to belong to these generators. Concerning Young diagrams and Young tableaux we use the notations and conventions of \cite{Ful} (see also in \cite{D-VH2}).  Concerning quantum groups, the quantum group $GL_q(2)$ is described for instance in \cite{YuM2} while the two parameter $GL_{p,q}(2)$ can be found in \cite{malt}. To make contact with the conventions of some authors one should replace $p$ by $p^{-1}$; here the convention is such that $GL_q(2)=GL_{q,q}(2)$. By an action of $GL_q(2)$ for instance we here mean a coaction of the Hopf algebra which is the appropriate deformation of the Hopf algebra of (representative) functions on $GL(2)$, i.e. the quantum group is a dual object to the corresponding Hopf algebra.

 \section{Homogeneous algebras, Koszulity}
 \setcounter{equation}{0}
 
 In this section we recall the definitions and results of \cite{BD-VW} and \cite{RB3} needed for this paper and we establish a formula involving the Poincar\'e series of a homogeneous algebra and the Euler characteristic of the corresponding Koszul complex which generalizes the one known for a quadratic algebra. When the homogeneous algebra is a Koszul algebra this formula simplifies and gives a criterion of koszulity which will be used in the next sections.\\

Throughout the following $N$ is an integer with $N\geq 2$, $\mathbb K$ is a commutative field and all the vector spaces and algebras are over the field $\mathbb K$.\\

A {\sl homogeneous algebra of degree $N$ or $N$-homogeneous algebra} is an algebra of the form \cite{BD-VW}.
\begin{equation}
\cala=A(E,R) = T(E)/(R)
\label{eq2.1}
\end{equation}
where $E$ is a finite-dimensional vector space, $T(E)$ is the tensor algebra of $E$ and $(R)$ is the two-sided ideal of $T(E)$ generated by a vector subspace $R$ of $E^{\otimes^N}$. When $N=2$ or $N=3$ we shall speak of quadratic or cubic algebra respectively. In view of the homogeneity of $R$, one sees that, by giving the degree 1 to the elements of $E$, $\cala$ is a graded algebra $\cala=\oplus_{n\in {\mathbb N}} \cala_n$ which is connected ($\cala_0=\mathbb K$) generated in degree 1 and such that the $\cala_n$ are finite-dimensional vector spaces. So the {\sl Poincar\'e series of} $\cala$
\begin{equation}
P_\cala(t) = \sum_n \dim(\cala_n)t^n
\label{eq2.2}
\end{equation}
is well defined for such a $N$-homogeneous algebra.\\

Given a $N$-homogeneous algebra $\cala=A(E,R)$, {\sl its dual} $\cala^!$ is defined to be \cite{BD-VW} the $N$-homogeneous algebra $\cala^!=A(E^\ast,R^\perp)$ where $E^\ast$ is the dual vector space of $E$ and where $R^\perp\subset E^{\ast\otimes^{N}}$ is the annihilator  of $R$, 
\[
R^\perp=\{ \omega\in (E^{\otimes^N})^\ast\>\> \vert \>\> \omega(x)=0,\>\> \forall x\in R\}
\]
with the canonical identification $E^{\ast\otimes^{N}}=(E^{\otimes^N})^\ast$. One has $(\cala^!)^!=\cala$.\\

As explained in \cite{BD-VW} to $N$-homogeneous algebras are canonically associated $N$-complexes which generalize the Koszul complexes of quadratic algebras. Let us recall the construction of the $N$-complex $K(\cala)$ associated with the $N$-homogeneous algebra $\cala$. One sets $K(\cala)=\oplus_n K_n(\cala)$ where the $K_n(\cala)$ are the left $\cala$-modules
\begin{equation}
K_n(\cala)=\cala\otimes (\cala^!_n)^\ast
\label{eq2.3}
\end{equation}
for $n\in \mathbb N$. One has $(\cala^!_n)^\ast=E^{\otimes^n}$ for $n<N$ and 
\begin{equation}
(\cala^!_n)^\ast=\cap_{r+s=n-N} E^{\otimes^r}\otimes R\otimes E^{\otimes^s}
\label{eq2.4}
\end{equation}
for $n\geq N$. Thus one always has $(\cala^!_n)^\ast\subset E^{\otimes^n}$ and the (left) $\cala$-module homomorphisms of $\cala\otimes E^{\otimes^{n+1}}$ into $\cala\otimes E^{\otimes^n}$ defined by 
\begin{equation}
a\otimes (e_0\otimes e_1\otimes \dots \otimes e_n)\mapsto (ae_0)\otimes (e_1\otimes \dots \otimes e_n)
\label{eq2.5}
\end{equation}
induce $\cala$-module homomorphisms $d:K_{n+1}(\cala)\rightarrow K_n(\cala)$ for any $n\in \mathbb N$. One has $d^N=0$ on $K(\cala)$ so $K(\cala)$ is a $N$-complex of left $\cala$-modules. One has
\begin{equation}
d(\cala_r\otimes (\cala_{s+1}^!)^\ast)\subset \cala_{r+1}\otimes (\cala^!_s)^\ast
\label{eq2.6}
\end{equation}
so $K(\cala)$ splits into sub $N$-complexes
\begin{equation}
K^{(n)}(\cala)=\oplus_m \cala_{n-m} \otimes (\cala^!_m)^\ast
\label{eq2.7}
\end{equation}
which are homogeneous for the total degree. \\

From the chain $N$-complex of left $\cala$-modules $K(\cala)$ one obtains by duality the cochain $N$-complex of right $\cala$-modules $L(\cala)$, see in \cite{BD-VW}.\\

In the case $N=2$ that is when $\cala$ is a quadratic algebra $K(\cala)$ is a complex, {\sl the Koszul complex of} $\cala$, the acyclicity of which in  positive degrees characterizes the quadratic {\sl Koszul algebras}. This is a class of very regular algebras which contains the algebras of polynomials. It is natural to look for a generalization of the notion of koszulity for $N$-homogeneous algebras with arbitrary $N\geq 2$. As shown in \cite{BD-VW} the complete acyclicity in positive degrees of the $N$-complex $K(\cala)$ is too strong for $N\geq 3$. In fact a good generalization of the notion of koszulity for a $N$-homogeneous  algebra $\cala$ was defined in \cite{RB3} and characterized there by the acyclicity in positive degrees of an ordinary complex which was identified in \cite{BD-VW} with one of the complexes obtained by contraction of the $N$-complex $K(\cala)$ (i.e. by setting alternatively $d^p,d^{N-p}$ for the differential). Furthermore it was shown in \cite{BD-VW} that (for $N\geq 3$) this complex is the only complex obtained by contraction of $K(\cala)$ for which the acyclicity in positive degrees leads to a non trivial class of algebras. This complex, which following \cite{RB3} will be refered to {\sl as the Koszul complex of $\cala$}, is obtained by putting together alternatively $N-1$ or $1$ arrows $d$ of $K(\cala)$ and by starting as
\begin{equation}
\dots \stackrel{d}{\longrightarrow} \cala\otimes (\cala^!_N)^\ast \stackrel{d^{N-1}}{\longrightarrow} \cala\otimes (\cala^!_1)^\ast \stackrel{d}{\longrightarrow} \cala \longrightarrow 0
\label{eq2.8}
\end{equation}
It is denoted by $C_{N-1,0}=C_{N-1,0}(K(\cala))$. The (other) contraction of $K(\cala)$ given by
\[
\dots \stackrel{d^{N-p}}{\longrightarrow} \cala\otimes (\cala^!_{N+r})^\ast \stackrel{d^p}{\longrightarrow} \cala\otimes (\cala^!_{N-p+r})^\ast \stackrel{d^{N-p}}{\longrightarrow} \cala\otimes (\cala^!_r)^\ast \longrightarrow 0
\]
for $0\leq r<p\leq N-1$ is denoted by $C_{p,r}=C_{p,r}(K(\cala))$.\\

The splitting of the $N$-complex $K(\cala)$ induces a corresponding splitting of the complex $C_{N-1,0}$ into subcomplexes $C^{(n)}_{N-1,0}$
\[
C_{N-1,0}=\oplus_n C^{(n)}_{N-1,0}
\]
which are homogeneous for the total degree with $C^{(n)}_{N-1,0}$ given by
\begin{equation}
\dots \stackrel{d}{\longrightarrow} \cala_{n-N} \otimes (\cala^!_N)^\ast \stackrel{d^{N-1}}{\longrightarrow} \cala_{n-1} \otimes (\cala^!_1)^\ast \stackrel{d}{\longrightarrow} \cala_n\longrightarrow 0
\label{eq2.9}
\end{equation}
with obvious conventions (e.g. $\cala_k=0$ for $k<0$). The Euler characteristic $\chi^{(n)}$ of $C^{(n)}_{N-1,0}$ can be computed in terms of the dimensions of the $\cala_r$ and $\cala^!_s$
\begin{equation}
\chi^{(n)}=\sum_{k\geq 0} (\dim (\cala_{n-kN})\dim (\cala^!_{kN})-\dim (\cala_{n-kN-1})\dim (\cala^!_{kN+1}))
\label{eq2.10}
\end{equation}
and it is worth noticing that the $\chi^{(n)}$ are finite; in fact the complexes  $C^{(n)}_{N-1,0}$ are finite-dimensional since they are of finite length and each term $\cala_k\otimes (\cala^!_\ell)^\ast$ is finite-dimensional. Let us define the series
\begin{equation}
\chi_\cala(t)=\sum_n \chi^{(n)} t^n
\label{eq2.11}
\end{equation}
\begin{equation}
Q_\cala(t)=\sum_n(\dim(\cala^!_{nN})t^{nN}-\dim(\cala^!_{nN+1})t^{nN+1})
\label{eq2.12}
\end{equation}
in terms of which one has the following result in view of (\ref{eq2.10}).

\begin{proposition}
One has the following relations for a $N$-homogeneous algebra $\cala$
\[
P_\cala(t) Q_\cala(t)=\chi_\cala(t).
\]
\end{proposition}

By definition, a $N$-homogeneous algebra is {\sl Koszul} \cite{RB3} if the Koszul complex $C_{N-1,0}(\cala)$ is acyclic in positive degrees or which is the same if the complexes $C^{(n)}_{N-1,0}$ are acyclic for $n>0$ (i.e. $n\geq 1$). It is clear that $\chi^0=1$ so if $\cala$ is Koszul $\chi^{(n)}=0$ for $n\geq 1$ and one has $\chi_\cala(t)=1$ which implies that Proposition 1 has the following corollary.

\begin{corollary} 
 
Let $\cala$ be a $N$-homogeneous algebra which is Koszul then one has
\[
P_\cala(t) Q_\cala(t)=1.
\]
\end{corollary}

It is worth noticing here that this corollary can also be deduced from Proposition 2.9 of \cite{ATVB2}.
In the case $N=2$, i.e. when $\cala$ is quadratic, one has $Q_\cala(t)=P_{\cala^!}(-t)$ so one sees that Proposition 1 and Corollary 1 above generalize well known results for quadratic algebras \cite{Pri}, \cite{YuM2}.\\

Corollary 1 is very useful to compute the Poincar\'e series $P_\cala(t)$ of a Koszul algebra $\cala$ when $\cala^!$ is small and easy to compute (which is frequent), it is also useful in order to prove that a $N$-homogeneous algebra $\cala$ is not Koszul when one knows $P_\cala(t)$ and $Q_\cala(t)$. It is this latter application that will be used in the next sections while the former application is used in \cite{AC.MDV2}. 

\section{The parafermionic algebra $\calb$}
\setcounter{equation}{0}

From now on we specialize to the case where $\mathbb K$ is the field $\mathbb C$ of complex numbers. By using the notations of Section 2, one can define $\calb$ to be the cubic algebra
\begin{equation}
\calb = A(\mathbb C^D, R_\calb)
\label{eq3.1}
\end{equation}
with $R_\calb\subset (\mathbb C^D)^{\otimes^3}$ defined to be the linear span of the set 
\begin{equation}
\{\left[ [x,y]_\otimes ,z\right]_\otimes\>\>  \vert\>\>  x,y,z \in \mathbb C^D\}
\label{eq3.2}
\end{equation}
where $[x,y]_\otimes=x\otimes y - y\otimes x$. It is clear that $R_\calb$ is invariant by the action of the linear group $GL(D)=GL(D,\mathbb C)$. This implies that $GL(D)$ acts on the algebra $\calb$ by automorphisms which preserve the degree. The irreducible representations of $GL(D)$ occuring here are labelled by
 the Young diagrams (i.e. by the partitions) and each one appears with multiplicity 1, \cite{OK} (see also in \cite{Ful} Exercise 15). In the representation space $\calb^\lambda$ corresponding to the Young diagram $\lambda$, there is a linear basis labelled by the Young tableaux of shape $\lambda$ and therefore $\calb(=\oplus_\lambda \calb^\lambda)$ admits as announced in the introduction a homogeneous basis labelled by the Young tableaux; in degree $n$ one has $\calb_n=\oplus_{\vert \lambda\vert=n} \calb^\lambda$ where $\vert\lambda\vert$ denotes the number of cells of the Young diagram $\lambda$. Using this, one can compute the dimensions of the homogeneous components. The most compact and useful form is the resulting Poincar\'e series of $\calb$ which is given by \cite{chat}
\begin{equation}
P_\calb(t) = \left(\frac{1}{1-t}\right)^D \left(\frac{1}{1-t^2}\right)^{\frac{D(D-1)}{2}}
\label{eq3.3}
\end{equation} 
From Formula (\ref{eq3.3}) one deduces that $\calb$ has polynomial growth. In fact one has
\[
\dim (\calb_n)\sim Cn^{\frac{D(D+1)}{2}-1}
\]
that is $gk$-$\dim(\calb)=\frac{D(D+1)}{2}$.\\

\noindent \underbar{Remark}. Another instructive way to establish Formula (\ref{eq3.3}) is the following. Consider $\mathbb C^D\oplus \wedge^2 \mathbb C^D$ equipped with the bracket defined by $[x,y]=x\wedge y$ if $x$ and $y$ are both in $\mathbb C^D$ and $[x,y]=0$ otherwise. This is a Lie bracket and $\mathbb C^D\oplus \wedge^2 \mathbb C^D$ is a graded Lie algebra for this bracket if one gives the degree 1 to the elements of $\mathbb C^D$ and the degree 2 to the elements of $\wedge^2\mathbb C^D$. By definition $\calb$ is the universal enveloping algebra of this graded Lie algebra and it is isomorphic as graded vector space to the symmetric algebra 
\[
S(\mathbb C^D\oplus \wedge^2\mathbb C^D)=S(\mathbb C^D)\otimes S(\wedge^2\mathbb C^D)
\]
with graduation induced by the one of $\mathbb C^D \oplus \wedge^2 \mathbb C^D$. Formula (\ref{eq3.3}) follows immediately.\\

Let us now come to the description of the dual cubic algebra $\calb^!=A(\mathbb C^{D\ast},R^\perp_\calb)$ of $\calb$. It is not hard to verify that $R^\perp_\calb\subset (\mathbb C^{D\ast})^{\otimes^3}$ is the linear span of the set
\begin{equation}
\{ \alpha\otimes \beta\otimes\gamma - \gamma\otimes \beta \otimes \alpha,\> \theta^{\otimes^3} \vert \alpha, \beta, \gamma, \theta \in \mathbb C^{D\ast}\}
\label{eq3.4}
\end{equation}
by using the obvious fact that $\dim(R_\calb)+\dim(R^\perp_\calb)=D^3$. Thus $\calb^!$ is generated (in degree 1) by the elements of $\mathbb C^{D\ast}=\calb^!_1$ with relations
\begin{equation}
\alpha\beta\gamma = \gamma\beta\alpha\>\>\> \mbox{and}\>\>\> \theta^3=0,\>\> \forall\alpha,\beta, \gamma, \theta \in \mathbb C^{D\ast}
\label{eq3.5}
\end{equation}
which implies that the symmetrized product and the antisymmerized product of 3 elements of $\calb^!_1=\mathbb C^{D\ast}$ vanish and that the product of 5 elements of $\calb^!_1=\mathbb C^{D\ast}$ also vanishes. Therefore one has $\calb^!_n=0$ for $n\geq 5$. In fact $GL(D)$ also acts on $\calb^!$ by automorphisms which preserve the degree and the content in irreducible subspaces is given by
\begin{equation}
(\bullet)_0 \oplus \left(\,
\begin{tabular}{|c|}\hline
\\ \hline
 \end{tabular}\, \right)_1 \oplus 
 \left (\, 
 \begin{tabular}{|c|}\hline
\\ \hline
\\ \hline
 \end{tabular} \oplus \begin{tabular}{|c|c|}\hline
 & \\ \hline
 \end{tabular}\, \right )_2 \oplus 
 \left( \,
 \begin{tabular}{|c|c|}\hline
 & \\ \hline
 \\
 \cline{1-1}
 \end{tabular}
 \, \right )_3 \oplus 
 \left(\, 
 \begin{tabular}{|c|c|}\hline
 & \\ \hline
 & \\ 
 \hline
 \end{tabular}
 \, \right)_4
 \label{eq3.6}
  \end{equation}
  where the parenthesis corresponds to the homogeneous component and where $\bullet$ is the empty Young diagram corresponding to the trivial 1-dimensional representation. It follows that the Poincar\'e series of $\calb^!$ is given by 
\[
P_{\calb^!}(t) =1+Dt +D^2 t^2 + \frac{1}{3} D(D^2-1)t^3 
+\frac{1}{12}D^2(D^2-1)t^4
\]
while $Q_\calb(t)$ is given by
\begin{equation}
Q_\calb(t) =1-Dt+\frac{1}{3} D(D^2-1)t^3-\frac{1}{12}D^2(D^2-1) t^4
\label{eq3.7}
\end{equation}
so that by applying Proposition 1 one obtains 
\begin{equation}
\chi_\calb (t)=\frac{1-Dt+\frac{1}{3}D(D^2-1)t^3-\frac{1}{12}D^2(D^2-1) t^4}{(1-t)^D (1-t^2)^{\frac{D(D-1)}{2}}}
\label{eq3.8}
\end{equation}
for $\chi_\calb (t)$.\\

For $D=2$, one has $\chi_\calb(t)=1$ but in this case $\calb$ is a well known cubic Artin-Schelter regular algebra \cite{AS}, \cite{ATVB} which is therefore a Koszul algebra  \cite{RB3} of global dimension 3 and is Gorenstein (definition e.g. in Appendix of \cite{AC.MDV2}). In fact it is the universal enveloping algebra of the Heisenberg Lie algebra. For $D\geq 3$ then $\chi_\calb(t)\not= 1$ so $\calb$ cannot be a Koszul algebra in view of Corollary 1. It is however worth noticing that even in these cases, the Euler characteristics $\chi^{(1)}, \chi^{(2)}, \chi^{(3)}$ and $\chi^{(4)}$ vanish. In fact one has
\begin{equation}
\chi_\calb(t) = 1+ \frac{1}{20} D(D^2-1)(D^2-4)t^5 \Lambda_\calb (t)
\label{eq3.9}
\end{equation}
where the series $\Lambda_\calb(t)$ satisfies $\Lambda_\calb(0)=1$ for $D\geq 3$. On the other hand we have seen above that $\calb$ has polynomial growth for any $D$.

\section{The plactic algebra $\calp$}
\setcounter{equation}{0}

In the following we let $(e_k)$, $k\in \{1,\dots,D\}$ be the canonical basis of $\mathbb C^D$. One can define the algebra $\calp$ to be the cubic algebra
\begin{equation}
\calp = A(\mathbb C^D,R_\calp)
\label{eq4.1}
\end{equation}
with $R_\calp\subset (\mathbb C^D)^{\otimes^3}$ defined to be the linear span of
\[
\{e_\ell \otimes e_m \otimes e_k-e_\ell \otimes e_k \otimes e_m \vert k<\ell \leq m\} \cup \{e_k\otimes e_m \otimes e_\ell-e_m \otimes e_k \otimes e_\ell \vert k\leq \ell <m\}.
\]
In contrast to $R_\calb$, $R_\calp$ depends on the basis $(e_k)$ and even on the ordered set $\{1,\dots,D\}$. Thus there is no natural action of $GL(D)$ on $\calp$. Nevertheless $\calp$ is the algebra (over $\mathbb C$) of an associative monoid the elements of which are the Young tableaux \cite{Ful}. It follows that it admits like $\calb$ a homogeneous linear basis labelled by the Young tableaux. Therefore it has in particular the same Poincar\'e series as $\calb$, i.e. $P_\calp(t)=P_\calb(t)$. Thus $\calp$ has polynomial growth and the same $gk$-dimension as $\calb$.\\

We now describe the dual cubic algebra $\calp^!$ of $\calp$. Let $(\theta^k)$, $k\in \{1,\dots,D\}$ be the basis of $\mathbb C^{D\ast}$ dual to the basis $(e_\ell)$ of $\mathbb C^D$, i.e. such that $\langle \theta^k, e_\ell\rangle=\delta_{k\ell}$. One has
\[
\calp^!=A\left (\mathbb C^{D\ast},R^\perp_\calp\right)
\]
and $R^\perp_\calp \subset (\mathbb C^{D\ast})^{\otimes^3}$ is spanned by the elements
\begin{equation}
\left.
\begin{array}{lll}
\theta^j\otimes \theta^k\otimes \theta^i+\theta^j\otimes \theta^i \otimes \theta^k & \mbox{with} & i<j\leq k\\
\theta^i\otimes \theta^k \otimes \theta^j + \theta^k \otimes \theta^i \otimes \theta^j & \mbox{with} & i\leq j < k\\
\theta^i\otimes \theta^j \otimes \theta^k & \mbox{with} & i\leq j\leq k\\
\theta^k \otimes \theta^j \otimes \theta^i & \mbox{with} & i<j<k
\end{array}
\right\}
\label{eq4.2}
\end{equation}
Using (\ref{eq4.2}) one sees that $\calp^!_n=0$ for $n\geq 5$ and that moreover, in an obvious sense, $\calp^!$ has the same content (\ref{eq3.6}) in Young diagrams as $\calb^!$ (i.e. homogeneous linear basis labelled by the corresponding Young tableaux). So one also has
\begin{eqnarray}
P_{\calp^!}(t) & = & P_{\calb^!}(t)\\
Q_{\calp}(t) & = & Q_{\calb}(t)\\
\chi_{\calp}(t) & = & \chi_{\calb}(t)
\end{eqnarray}
and one can again conclude that for $D\geq 3$ the cubic algebra $\calp$ is not a Koszul algebra.\\

For $D=2$ a sharp difference appears between the cubic algebras $\calp$ and $\calb$ in that $\calp$ is not a cubic Artin-Schelter regular algebra. Let us investigate more closely this case. By definition the sequence
\begin{equation}
0\rightarrow \calp \otimes \calp^{!\ast}_4\stackrel{d}{\rightarrow} \calp \otimes R_\calp \stackrel{d^2}{\rightarrow}\calp\otimes \calp_1 \stackrel{d}{\rightarrow}\calp\rightarrow \mathbb C \rightarrow 0
\label{eq4.6}
\end{equation}
is exact at each term on the right of $\calp\otimes R_\calp$ for any value of $D$. On the other hand for $D=2$ one has $\chi_\calp(t)=1$ thus (\ref{eq4.6}) is exact for $D=2$ if and only if $d:\calp\otimes \calp^{!\ast}_4\rightarrow \calp\otimes R_\calp$ is injective which follows from the fact that for $D=2$, (\ref{eq4.6}) reads
\begin{equation}
0\rightarrow \calp \stackrel{
\begin{array}{c}
(0,e_2)\\
\end{array}}{\longrightarrow}\calp^2 \stackrel{
\left(\begin{array}{cc}
e^2_2 & -e_2e_1\\
\left[e_1,e_2\right] & 0
\end{array}\right )}{\hbox to 12mm{\rightarrowfill}}\calp^2 \stackrel{
\left(
\begin{array}{c}
e_1\\
e_2
\end{array}
\right )}{\longrightarrow} \calp \rightarrow \mathbb C\rightarrow 0
\label{eq4.7}
\end{equation}
the exactness of which is anyhow straightforward. It follows that for $D=2$ the cubic algebra $\calp$ is Koszul of global dimension 3 but the asymmetry of (\ref{eq4.7}) shows that it cannot be Gorenstein and thus that it is not regular\cite{AS}. As will be shown in \cite{TP}, this non Gorenstein property of $\calp$ for $D=2$ manifests itself by the non triviality of the cohomology of a cochain complex obtained by contraction of the cochain 3-complex of right $\calp$-modules $L(\calp)$ which is dual of $K(\calp)$  (see in \cite{BD-VW}).

\section{Parabosonic algebra and further prospect}
\setcounter{equation}{0}

As explained in the introduction we have restricted attention in Section 3 to the parafermionic case. The analysis in the parabosonic case is completely similar. Let us say a word on it. In the parabosonic case the analog of $\calb$ is the cubic algebra $\tilde \calb=A(\mathbb C^D, R_{\tilde\calb})$ with $R_{\tilde\calb}\subset (\mathbb C^D)^{\otimes^3}$ defined to be the linear span of the set
\begin{equation}
\{[\{x,y\}_\otimes, z]_\otimes\>\> \vert\>\>  x,y,z\in \mathbb C^D\}
\label{eq5.1}
\end{equation}
where $\{x,y\}_\otimes =x\otimes y + y\otimes x$. The content of $\tilde \calb$ in terms of representations (up to isomorphisms) of $GL(D)$ is the same as the one of $\calb$ so one can analyse it in the same manner. However it is more interesting to proceed here in the line of the remark of Section 2. In fact $\tilde\calb$ is the ``super" version of $\calb$. More precisely, consider $\mathbb C^D \oplus S^2 \mathbb C^D$ equipped with the (super) bracket defined by $[\![x,y]\!] =x\vee y$ if $x$ and $y$ are both in $\mathbb C^D$ and $[\![x,y]\!] =0$ otherwise. This is a graded super Lie algebra and by its very definition $\tilde\calb$ is its universal enveloping algebra. The super version of the Poincar\'e-Birkhoff-Witt theorem (see e.g. in \cite{Kost}) implies that as graded vector space (and coalgebra) $\tilde\calb$ is isomorphic to the ``super" symmetric algebra
\[
\tilde S (\mathbb C^D \oplus S^2\mathbb C^D)=\wedge (\mathbb C^D)\otimes S(S^2\mathbb C^D)
\]
with graduation induced by the one of $\mathbb C^D\oplus S^2\mathbb C^D$. It follows that the Poincar\'e series of $\tilde\calb$ is given by
\begin{equation}
P_{\tilde\calb}(t)=(1+t)^D \left(\frac{1}{1-t^2}\right)^{\frac{D(D+1)}{2}}
\label{eq5.2}
\end{equation}
from which
\begin{equation}
P_{\tilde \calb}(t) = P_\calb(t)
\label{eq5.3}
\end{equation}
follows immediately in view of (\ref{eq3.3}). The structure of $\tilde \calb^!$ is similar to the one of $\calb^!$ in particular they have the same Poincar\'e series (polynomial here) and one has again the relation
\begin{equation}
\chi_{\tilde\calb}(t)=\chi_\calb(t)
\label{eq5.4}
\end{equation}
connecting their Koszul homologies.\\

The above coincidences are not accidental, the algebras $\calb$ and $\tilde \calb$ belong to a family of algebras which depend continuously on parameters while $\calp$  is a limit in this family or, more precisely, is obtained for a singular value of these parameters. The complete general analysis (as well as the analysis of the invariance properties) of this family of algebras will be reported elsewhere \cite{TP}. In the following we just give the description of the family in the case $D=2$ where it relies to the Artin-Schelter classification \cite{AS}.\\

The generic form of the relations for the cubic Artin-Schelter regular algebra of type $S_1$ is (cf. (8.6) in \cite{AS})
\begin{equation}
\left.
\begin{array}{l}
e_2 e^2_1 + qr e^2_1 e_2 - (q+r) e_1 e_2 e_1 =0\\
\\
e^2_2 e_1 +qr e_1 e^2_2 - (q+r) e_2 e_1 e_2 =0
\end{array}
\label{eq5.5}
\right\}
\end{equation}
where $q$ and $r$ are complex parameters with $qr\not= 0$. The contact with the notations of \cite{AS} is obtained by the substitutions $qr\mapsto \alpha$, $-(q+r)\mapsto a$,  $e_1\mapsto x$ and $e_2\mapsto y$ which means in particular that $q$ and $r$ are the roots of the equation $t^2+at+\alpha=0$. We denote by $\cala_{q,r}$ the algebra generated $e_1$ and $e_2$ with relations (\ref{eq5.5}) above. For $qr\not=0$ this is a generic cubic Artin-Schelter regular algebra of type $S_1$ while for $qr=0$ this is a Kozsul algebra of global dimension 3 which is not Gorenstein \cite{AS}, (the argument being the same as in last section). One has $\calb=\cala_{1,1}$, $\tilde \calb=\cala_{-1,1}$ while $\calp=\cala_{0,1}$ is a singular case $qr=0$. Setting $r=1$, one obtains the 1-parameter subfamily $\cala_q=\cala_{q,1}$ which already contains $\calb$ and $\tilde \calb$ as well as $\calp$ which is the singular point $q=0$. In this case $r=1$, relations (\ref{eq5.5}) read 
\begin{equation}
[e_1e_2-q^{-1}e_2e_1,e_k]=0,\>\>\> k\in \{1,2\}
\label{eq5.6}
\end{equation}
which just expresses the fact that $e_1e_2-q^{-1}e_2e_1$ is central. In this form, their $GL_q(2)$ invariance is obvious for $q\not=0$ so one sees that the $GL(2)$ symmetry of $\calb$ for $D=2$ generalizes as symmetry of $\cala_q$ by the quantum group $GL_q(2)$ for $q\not=0$. This of course does not make sense for $q=0$ i.e. for $\calp$.  In fact it is easily verified that this extends for $qr\not=0$ as symmetry of $\cala_{q,r}$ by the quantum group $GL_{p,q}(2)$ with $p=q/r^2$ (with these notations one has $GL_q(2)=GL_{q,q}(2)$).It is worth noticing here that by definition one has $\cala_{q,r}=\cala_{r,q}$ so together with the quantum symmetry $GL_{q/r^2,q}(2)$ one also has the quantum symmetry $GL_{r/q^2,r}(2)$ for $qr\not=0$.\\

As will be explained in \cite{TP} where more details will be given,  this analysis of (quantum) symmetry generalizes in arbitrary dimension $D$ for the appropriate family of algebras containing $\calb$ and $\tilde\calb$ and admitting $\calp$ as limit or singular point.

\section*{Acknowledgements}

It is a pleasure to thank Roland Berger, Alain Connes, Oleg Ogievetsky and Ivan Todorov for stimulating discussions. 
 
\baselineskip=0,5cm

\end{document}